\documentclass[11pt,reqno]{amsart}

\usepackage{a4,graphicx,mathrsfs,amssymb, tikz}
\usepackage[colorlinks=true]{hyperref}
\usepackage{enumerate}

\newtheorem{theorem}{Theorem}[section]
\newtheorem{lemma}[theorem]{Lemma}

\newtheorem{proposition}[theorem]{Proposition}
\theoremstyle{definition}

\newtheorem{definition}[theorem]{Definition}
\theoremstyle{remark}

\numberwithin{equation}{section}

\renewcommand{\mathcal}{\mathscr}

\renewcommand{\hat}{\widehat}

\allowdisplaybreaks
\begin{document}

\title[Exponential stability of a transmission problem]{Exponential stability of a damped beam-string-beam transmission problem}\thanks{The authors would like to thank
COLCIENCIAS (Project 121571250194) for the financial support.}

\author[B. Barraza]{Bienvenido Barraza Mart\'inez}
\address{B.\ Barraza Mart\'inez, Universidad del Norte, Departamento de Matem\'aticas y Estad\'istica, Barranquilla, Colombia}
\email{bbarraza@uninorte.edu.co}

\author[J. Hern\'andez]{Jairo Hern\'andez Monz\'on}
\address{J.\ Hern\'andez Monz\'on, Universidad del Norte, Departamento de Matem\'aticas y Estad\'istica, Barranquilla, Colombia}
\email{jahernan@uninorte.edu.co}

\author[G. Vergara]{Gustavo Vergara Rolong}
\address{G.\ Vergara Rolong, Universidad de la Costa, Departamento de Ciencias Naturales y Exactas, Barranquilla, Colombia}
\email{gvergara@cuc.edu.co}

\date{\today}

\begin{abstract}
We consider a beam-string-beam transmission problem, where two structurally
damped or undamped beams are coupled with a frictionally damped string by
transmission conditions. We show that for this type of structure, the
dissipation produced by the frictional part is strong enough to produce
exponential decay of the solution, no matter how small is its size: for the
exponential stability in the damped-damped-damped situation we use energy
method and in the undamped-damped-undamped situation we use a frequency domain
method from the semigroups theory, which combines a contradiction argument
with the multiplier technique to carry out a special analysis for the
resolvent. Additionally, we show that the solution first defined by the weak
formulation, in fact have higher Sobolev space regularity.
\end{abstract}


\maketitle

\section{Introduction}
\setcounter{equation}{0}
\bigskip Recent advances in material science have provided new means of
suppressing vibrations from elastic multi-link structures, for instances, by
applying some type of local or total dampings to them. These structures
consisting of connected flexible elements such as strings, beams, plates and
shells have many applications in engineering areas such as in robot arms,
frames, solar panels, aircrafts, satellite antennae, bridges and so on (see
\cite{BSW95}, \cite{BSW96}, \cite{LLS94} and the references therein). In this
context we consider, in this paper, a coupled beam-string-beam system, where
we assume structural damping/no damping for the beams and frictional damping
for the string. More precisely, we consider a elastic structure composed by
three parts. The firt and the third are structurally damped or undamped beams
in the open intervals $\left(  l_{0},l_{1}\right)  $ and $\left(  l_{2}%
,l_{3}\right)  $, respectively, and the second is a frictionally damped
string, occupying in equilibrium the open interval $\left(  l_{1}%
,l_{2}\right)  $, as shown in Figure 1.\\
\begin{center}
    \includegraphics[scale=1.5]{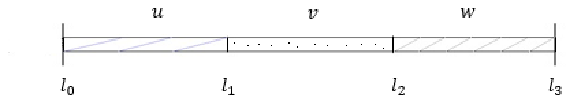}\\
    Figure 1
\end{center}
Denoting by $u=u(x,t)$, $w=w(x,t)$ and $v=v(x,t)$ the vertical displacements
of the points on the beams and on the string with coordinates $x$ at time $t$,
respectively, the mathematical model for the structure is given by the
equations\\
\begin{align}
     u_{tt}+u_{xxxx}-\rho _1 u_{txx}=0,\quad \text{in}\   (l_0 , l_1) \times (0,  \infty), \label{1}\\
    v_{tt}-v_{xx}+\beta v_{t}=0, \quad \text{in}\  (l_1 , l_2) \times (0,  \infty),\label{2}\\
    w_{tt}+w_{xxxx}-\rho _2 w_{txx}=0, \quad \text{in}\  (l_2 , l_3) \times (0,  \infty),\label{3}
\end{align}
where $\rho_{i}\geq0$, $i=1,2$, and $\beta\geq0$ are fixed constants. The
coefficients $\rho_{1}\geq0$ and $\rho_{2}\geq0$ describe the structural
damping (or the absence of damping) for the beam equations \eqref{1} and \eqref{3},
whereas $\beta>0$ in \eqref{2} describes a frictional damping on the string. On
the endpoints $l_{0}$, $l_{3}$, of the beams, we impose clamped (Dirichlet)
boundary conditions
\begin{align}
    u(l_0,t)=u_x(l_0,t)=w(l_3,t)=w_x (l_3,t)=0,\qquad     t \in (0,\infty).\label{4}
\end{align}   
On the interface $\left\{  l_{1},l_{2}\right\}  $, we have transmission
conditions
\begin{align}
    u(l_1,t)=v(l_1, t)\text{ and } v(l_2,t)=w(l_2,t), \qquad  t \in (0,\infty),\label{6}\\
    u_{xxx}(l_1,t)-\rho _1 u_{tx}(l_1,t)+v_x (l_1,t)=0,  \qquad t \in (0,\infty),\label{7} \\
     w_{xxx}(l_2,t)-\rho _2 w_{tx}(l_2,t)+v_x (l_2,t)=0, \qquad  t \in (0,\infty),\label{8}\\
    u_{xx}(l_1,t)=0,\qquad \qquad  t \in (0,\infty),\label{9}\\
    w_{xx}(l_2,t)=0, \qquad \qquad  t \in (0,\infty).\label{10}
\end{align}
Condition (\ref{6}) is known as the continuity transmission condition, (\ref{7}) and
(\ref{8}) mean that the two forces which are the shear force of the beams and the
stress of the string are such that one cancels the other, and (\ref{9}) and (\ref{10})
describe the fact that the the beams present possible inflection point on
$l_{1}$ and $l_{2}$ (compare with \cite{H16}, p. 1934).

Finally, the boundary-transmission problem (\ref{1})--(\ref{10}) is endowed with
initial conditions of the form
\begin{align}
     u(x,0)=u_0(x),\quad u_t(x,0)=u_1(x),\quad x\in (l_0 ,l _1), \\
     v(x,0)=v_0(x),\quad v_t(x,0)=v_1(x),\quad  x\in (l_1 ,l _2),\\
     w(x,0)=w_0 (x), \quad w_t(x,0)=w_1 (x),\quad  x\in (l_2 ,l _3).\label{5}
\end{align} 
\bigskip

The aim of the present paper is to study well-posedness, regularity, and
exponential stability of the solution of (\ref{1})-(\ref{5}).\\

In the recent years, the large-time behavior of structures consisting of
elastic strings and/or beams with different damping has been studied
a lot. We refer to \cite{LL98}, \cite{MP01} and \cite{Sh15}, where structures
formed by beams were studied. For instance, Shel in \cite{Sh15} showed, under
certain conditions, the exponential stability of a network of elastic and
thermoelastic Euler-Bernoulli beams. For  transmission problems between
strings, see for example, \cite{AMSV14}, \cite{BR10}, \cite{HZ17},
\cite{MMRV18}, \cite{RBA11} and \cite{RW2020}. For example, Alves, Mu\~{n}oz
Rivera et al. considered in \cite{AMSV14} a transmission problem of a material
composed of three components; one of them is a Kelvin--Voigt viscoelastic
material, the second is an elastic material (no dissipation), and the third is
an elastic material inserted with a frictional damping mechanism. They proved
exponential stability of the solution, if viscoelastic component is not in the
middle of the material. Then, Rissel and Wang established in \cite{RW2020}
exponential stability for a coupled system of elasticity and thermoelasticity
with second sound confined by a purely elastic one. Specifically, its
structure was composed of three parts, where the first and third were purely
elastic and the second thermo-elastic. Now, regarding elastic structures
composed of string and beam we mention \cite{AJM09}-\cite{AS18}, \cite{H16},
\cite{LHX18}, \cite{Sh20}, \cite{W19} and \cite{WW2020}. Ammari et al. in
\cite{AJM09}-\cite{AS18} considered the nodal feedback stabilization for
networks of strings and beams. They obtained that the decay rate of a
closed-loop system depends on the positions of the nodal feedback controllers.
Hassine in \cite{H16} studied a elastic transmission wave/beam systems with a
local Kelvin--Voigt damping. He showed that the energy of this coupled system
decays polynomially as the time variable goes to infinity, if the damping,
which is locally distributed, acts through one part
of the structure. Li, Han and Xu in \cite{LHX18} obtained polynomial
stability for a string-frictionally damped beam system and exponential
stability for a frictionally damped string-beam system. Shel in \cite{Sh20}
considered transmission problems for a coupling of a string and a beam with at
least one of them being thermoelastic and established that the associated
semigroup is exponentially stable when the string is thermoelastic and
polynomial stables when only the beam is thermoelastic and satisfies certain
additional condition. Wang in \cite{W19} obtained the strong stability of the
semigroup associated to a frictionally damped string-beam system. F. Wang and
J. M. Wang in \cite{WW2020} established exponential stability for a
beam-frictionally damped string system with some feedback at the interface point.\\

In this paper, we study the well-posedness of the problem (\ref{1})-(\ref{5}), the
higher regularity of the solution and the exponential stability of the energy
of the system, depending on the dampings. By energy method, we prove the
exponential stability of (\ref{1})-(\ref{5}), if the beams and the string are damped
(i.e. if $\rho_{1},\rho_{2}$ and $\beta$ are positive). For the proof we do
not need higher regularity of the solution. Moreover, by a frequency domain
method from the semigroups theory, we show the exponential stability of
(\ref{1})-(\ref{5}) in the undamped-damped-undamped situation, i.e. $\rho_{1}%
=\rho_{2}=0$ and $\beta>0$. For the proof we use the higher Sobolev space regularity of the solutions, which implies that the transmission conditions, first defined by the weak formulation, hold in the classical sense. For this we follow some analogous ideas from
Section 4 of \cite{BDHKN19}.\\

This paper is organized as follows: In Section 2, we define the basic
spaces and operators. Then, in Section 3 we show the generation of a $C_{0}%
-$semigroup of contractions (and therefore the well-posedness of
(\ref{1})-(\ref{5})). Exponential stability for the cases damped-damped-damped and
undamped-damped-undamped, as well as  higher regularity of the solutions, are
shown in Section 4.\\

Finally, let us say some words about notation. Derivatives with respect to $t$ of a function will be eventually denoted by a ``dot'' over the name of the fuction. So $\dot{\phi}$ will denote the derivative of $\phi$ with respect to $t$. As usual we will also use the notations $\psi'$, $\psi''$, $\psi'''$, or in general $\psi^{(n)}$ for $n\geq 4$, for the derivatives of a function $\psi$, which depends only on the one-dimensional spatial variable $x$.

\section{Base Spaces}
\setcounter{equation}{0}
For the existence and uniqueness of the solution of the problem  (\ref{1})-(\ref{5}) we will write this in an abstract form. For it, we will introduce some definitions and notations. For example, let us define the operator matrix $A$, which domain will be explained later,  by
 $$A:=\begin{pmatrix}
0 & 0 & 0 & 1 & 0 &0\\
0 & 0 & 0 & 0 & 1 &0\\
0 & 0 & 0 & 0 & 0 &1\\
-\dfrac{\partial ^4}{\partial x^4} & 0 & 0 & \rho _1\dfrac{\partial ^2}{\partial x^2} & 0 &0\\
0 & \dfrac{\partial ^2}{\partial x^2} & 0 & 0 & -\beta &0\\
0 & 0 & -\dfrac{\partial ^4}{\partial x^4} & 0 & 0 &\rho _2 \dfrac{\partial ^2}{\partial x^2}\\
 \end{pmatrix}.$$
 
 Then, if $z_t :=\dfrac{\partial z}{\partial t}$, $z_{tt} :=\dfrac{\partial ^2 z}{\partial t ^2}$ and $U:=(u,v,w,u_t,v_t,w_t)^\top$, we have that  $U_t=(u_t,v_t,w_t,u_{tt},v_{tt},w_{tt})^\top$ and (\ref{1})-(\ref{3}) can be written as
 $$U_t=\begin{pmatrix}
 u_{t}\\
 v_{t}\\
 w_{t}\\
u_{tt}\\
 v_{tt}\\
 w_{tt}\\
 \end{pmatrix}=\left(\begin{matrix}
 u_{t}\\
 v_{t}\\
 w_{t}\\
 -\small{\dfrac{\partial ^4 u}{\partial x^4} +\rho_1 \dfrac{\partial ^2 u_t}{\partial x^2}}\\[2ex]
  \small{\dfrac{\partial ^2 v}{\partial x^2}-\beta v_2 }\\[2ex]
  -\small{\dfrac{\partial ^4 w}{\partial x^4} +\rho_2 \dfrac{\partial ^2 w_t}{\partial x^2}}\\
 \end{matrix}\right)=AU .$$
So, with the initial condition (\ref{5}), we obtain
\begin{equation} \label{c}
   U_t(t)=AU(t)\quad (t>0), \qquad U(0)=U_0,
\end{equation}
where $U_0:=(u_0,v_0,w_0,u_1,v_1,w_1)^{\top}$.\\

Now, we define the following spaces
$$H_{l_0}^2 :=\{u\in 
H^2(I_1):u(l_0)=u'(l_0)=0\},$$
$$H_{l_3}^2 :=\{w\in 
H^2(I_3):w(l_3)=w'(l_3)=0\},$$
with the inner products
\begin{align}
    \langle u, \Tilde{u}\rangle_{H^2 _{l_0}} &:= \langle u'' , \Tilde{u}''\rangle _{L^2 (I_1)} \label{11},\\
      \langle w, \Tilde{w}\rangle_{H^2 _{l_3}} &:= \langle w'' , \Tilde{w}''\rangle _{L^2 (I_3)} \label{12}. \end{align}
Due to the generalized Poincar\'e's  inequality the induced norms $\Vert \cdot\Vert_{H^2_{l_0}}$ and  $\Vert\cdot\Vert_{H^2_{l_3}}$ are equivalents to the standard norms $\Vert\cdot\Vert_{H^2(I_1)}$ and $\Vert\cdot\Vert_{H^2(I_3)}$ on $H^2_{l_0}$ and $H^2_{l_3}$, respectively.\\

Now, let define the spaces
\begin{align*}
    \mathbb{H}&:=\{(u,v,w)^\top\in H_{l_0}^2 \times H^1(I_2) \times H_{l_3}^2: u(l_1)=v(l_1)\text{ y } v(l_2)=w(l_2) \},\\
    \mathbb{L}&:=L^2(I_1)\times L^2(I_2) \times L^2(I_3),
\end{align*}
equipped with the inner products
\begin{align}
    \langle (u,v,w)^\top,(\Tilde{u},\Tilde{v},\Tilde{w})^\top\rangle _{\mathbb{H}}&:= \langle u , \Tilde{u} \rangle _{H^2_{l_0}}+\langle v' , \Tilde{v}' \rangle _{L^2(I_2)}+ \langle w , \Tilde{w} \rangle _{H^2_{l_3}}, \label{p1}\\
    \langle (u,v,w)^\top,(\Tilde{u},\Tilde{v},\Tilde{w})^\top\rangle _{\mathbb{L}}&:= \langle u , \Tilde{u} \rangle _{L^2(I_1)}+\langle v , \Tilde{v} \rangle _{L^2(I_2)}+\langle w , \Tilde{w} \rangle _{L^2(I_3)}. \label{p2}
\end{align}\\
Again, by the Poincar\'e's inequality, the norm in $\mathbb{H}$, induced by the inner product \eqref{p1}, is equivalent to the standard norm in the poduct space $H^2(I_1)\times H^1(I_2)\times H^2(I_3)$. Due to the continuity of the trace operator, $\mathbb{H}$ is a closed subspace of $H^2(I_1)\times H^1(I_2)\times H^2(I_3)$ and therfore $\big(\mathbb{H},\langle\cdot , \cdot\rangle_{\mathbb{H}}\big)$ is a Hilbert space.\\

Finally, we define the Hilbert space
$$\mathcal{H}:=\mathbb{H}\times \mathbb{L},$$
with the inner product given by
\begin{align*}
    \langle U, \Tilde{U} \rangle _{\mathcal{H}}:&=\langle (u_1,v_1,w_1)^\top,(\Tilde{u}_1,\Tilde{v}_1,\Tilde{w}_1)^\top \rangle _{\mathbb{H}}+\langle (u_2,v_2,w_2)^\top,(\Tilde{u}_2,\Tilde{v}_2,\Tilde{w}_2)^\top \rangle _{\mathbb{L}} \\
    &=\langle u_1, \Tilde{u}_1 \rangle _{H^2 _{l_0}}+\langle v_{1}' , \Tilde{v}'_{1} \rangle _{L^2(I_2)}+\langle w_1 , \Tilde{w}_1 \rangle _{H^2 _{l_3}}+\langle u_2 , \Tilde{u}_2 \rangle _{L^2(I_1)}\\
    & \quad +\langle v_2 , \Tilde{v}_2 \rangle _{L^2(I_2)}+\langle w_2 , \Tilde{w}_2 \rangle _{L^2(I_3)}, 
\end{align*}
for all $U=(u_1,v_1,w_1,u_2,v_2,w_2)^{\top}$, $\Tilde{U}=(\Tilde{u}_1,\Tilde{v}_1,\Tilde{w}_1,\Tilde{u}_2,\Tilde{v}_2,\Tilde{w}_2)^{\top} \in \mathcal{H}$. \\

Since the functions in the spaces defined above are not regular enough to satisfy the transmission conditions (\ref{7})-(\ref{10}) in the classic sense, or even in the trace sense, we interpretate this transmission condition first in a ``weak'' sense. For this, we consider 
 $U=(u_1,v_1,w_1,u_2,v_2,w_2)^{\top} \in \mathcal{H}$ and $\Tilde{U}=(\Tilde{u}_1,\Tilde{v}_1,\Tilde{w}_1,\Tilde{u}_2,\Tilde{v}_2,\Tilde{w}_2)^{\top}\in \mathbb{H}\times \mathbb{H}$ sufficiently smooth, such that the following calculations make sense. Applying  integration by parts we obtain
\begin{align*}
    \langle AU, \Tilde{U} \rangle &= a(U,\Tilde{U})+b(U,(\Tilde{u}_2,\Tilde{v}_2,\Tilde{w}_2)^\top),
\end{align*}
where
\begin{align*}
    a(U,\Tilde{U})&:=\langle u_2 , \Tilde{u}_1 \rangle _{H^2 _{l_0}}+\langle v'_2 , \Tilde{v}'_1 \rangle _{L^2(I_2)}+\langle w_2 , \Tilde{w}_1 \rangle _{H^2 _{l_3}} - \langle u_1 '' \Tilde{u}''_2  \rangle _{L^{2}(I_1)}\\
    & \quad -\rho _1 \langle u_2 ', \Tilde{u}' _2 \rangle _{L^{2}(I_1)} -\beta \langle v_2 , \Tilde{v}_2 \rangle _{L^2(I_2)}
    -\langle v_1 ' , \Tilde{v}_2 '  \rangle _{L^{2}(I_2)} - \langle w_1 '', \Tilde{w}''_2  \rangle _{L^{2}(I_3)}\\
    & \quad -\rho _2 \langle w_2 ', \Tilde{w}' _2 \rangle _{L^{2}(I_3)}
\end{align*}
and
\begin{align*}
b(U,&(\Tilde{u}_2,\Tilde{v}_2,\Tilde{w}_2)^\top)\\
 & :=-u_{1}^{(3)}(l_1)\Tilde{u}_2 (l_1)+u_{1}^{(3)}(l_0)\Tilde{u}_2 (l_0) +u_{1}''(l_1)\Tilde{u}'_2 (l_1)-u_{1}''(l_0)\Tilde{u}'_2 (l_0)\\
& \quad +\rho _1 u' _2 (l_1) \Tilde{u}_2 (l_1)-\rho _1 u' _2 (l_0) \Tilde{u}_2 (l_0)+ v_1 ' (l_2)\Tilde{v}_2(l_2)-v_1 ' (l_1)\Tilde{v}_2(l_1)\\
    & \quad -w_{1}^{(3)}(l_3)\Tilde{w}_2 (l_3)+w_{1}^{(3)}(l_2)\Tilde{w}_2 (l_2)+w_{1}''(l_3)\Tilde{w}'_2 (l_3)-w_{1}''(l_2)\Tilde{w}'_2 (l_2)\\
    &\quad+\rho _2 w' _2 (l_3) \Tilde{w}_2 (l_3)-\rho _2 w' _2 (l_2) \Tilde{w}_2 (l_2).
\end{align*}
Since $\Tilde{U}\in \mathbb{H}\times \mathbb{H}$, we have that $ \Tilde{u}_2 (l_0)=\Tilde{u}' _2 (l_0)=\Tilde{w}_2 (l_3)=\Tilde{w}' _2 (l_3)=0$, $\Tilde{v}_2(l_1)=\Tilde{u}(l_1)$ y $\Tilde{v}_2(l_2)=\Tilde{w} _2(l_2)$. Then, it follows that
\begin{align*}
    b(&U,(\Tilde{u}_2,\Tilde{v}_2,\Tilde{w}_2)^\top)\\
    &=\big[-u_{1}^{(3)}(l_1)+\rho _1 u' _2 (l_1)-v_1 ' (l_1)\big]\Tilde{u}_2 (l_1)+\big[w_{1}^{(3)}(l_2)-\rho _2 w' _2 (l_2)+v_1 ' (l_2)\big]\Tilde{w}_2 (l_2)\\
    &\quad+u_1 '' (l_1)\Tilde{u}_2 ' (l_1)-w_1 '' (l_2)\Tilde{w}_2 ' (l_2)\\
    &=0,
\end{align*}
if and only if
\begin{align*}
    u_{1}^{(3)}(l_1)-\rho _1 u' _2 (l_1)+v_1 ' (l_1)&=0 ,\\
    w_{1}^{(3)}(l_2)-\rho _2 w' _2 (l_2)+v_1 ' (l_2)&=0,\\
    u_1 '' (l_1)&=0,\\
    w_1 '' (l_2)&=0,
\end{align*}
which are the transmission conditions (\ref{7})-(\ref{10}) that we have in the description of the problem in Section 1, if there $u_2=\partial_t u_1$ and $w_2=\partial_t w_1$. This motivates the following definition.
\vspace*{-4pt}
\begin{definition}
$U\in \mathcal{H}$ satisfies the transmission conditions (\ref{7})-(\ref{10}), in a weak sense, if and only if
\begin{align}
    \langle AU, \Tilde{U}  \rangle =a(U,\Tilde{U}), \text{ for all } \Tilde{U}\in \mathbb{H}\times \mathbb{H}.\label{cde}
\end{align}
\end{definition}
Now, we define the operator
\begin{align}
\mathcal{A}: &D(\mathcal{A})\subset \mathcal{H}\to \mathcal{H}, \quad U\to \mathcal{A}U:=AU,
\label{Ao}
\end{align}
where
\begin{align*}
    D(\mathcal{A}):=&\Big\{U\in \mathbb{H}\times \mathbb{H}:\text{ }u^{(4)}_ 1 \in L^2 (I_1), v''_1,\in L^2 (I_2), w^{(4)}_1\in L^2 (I_3), \\ &\quad U \text{ satisfies the (\ref{7})-(\ref{10}) conditions in a weak sense}   \Big\}.
\end{align*}
In this way, the problem (\ref{1})-(\ref{5}) can be written in abstract form by the Cauchy problem
\begin{equation}\label{cn}
    \dfrac{dU}{dt}(t)=\mathcal{A}U(t) \quad (t>0), \quad U(0)=U_0.
\end{equation}
\section{Well posedness }
\setcounter{equation}{0}
Now, we will show that the problem (\ref{cn}) is well posed, which means, that for every $U\in D(\mathcal{A})$, (\ref{cn}) has one and only classical solution which  continuously depends of $U_0$. For it, we will proof that the operator $\mathcal{A}$ defined on (\ref{Ao}) generates a  $C_0$-semigroup of contractions on $\mathcal{H}$. To achieve that, we will use the Lumer-Phillips theorem. 
\begin{proposition}\label{Proposition_wellposedness} The following assertions hold:
\begin{itemize}
    \item[a)] $\mathcal{A}$ is dissipative.
    \item[b)] $I- \mathcal{A}$ is surjective.
    \item[c)] $D(\mathcal{A})$ is dense in $\mathcal{H}$.
\end{itemize}
\end{proposition}
\begin{proof}
 Let $U\in D(\mathcal{A})$. From  (\ref{cde}) we have
\begin{align}\label{Eq_dissipativity_A}
    \langle \mathcal{A}U,U \rangle _{\mathcal{H}} 
    &=-\rho _1|| u_2' ||^{2}_{L^{2}(I_1)}-\beta ||v_2 ||^{2}_{L^{2}(I_2)}-\rho _2 || w_2'||^{2}_{L^{2}(I_3)} \leq 0,
\end{align}
from which $a)$ follows. For the surjectivity we will use the Lax-Milgram theorem.
Let $F=(f_1, g_1 , h_1 , f_2, g_2 , h_2)^{\top}\in \mathcal{H}$. We need to show that there exists a $U=(u_1, v_1 , w_1 , u_2, v_2 , w_2)^{\top}\in D(\mathcal{A})$ such that $ ( I -\mathcal{A})U=F $, i.e.
\begin{align}
    u_1 - u_2 &= f_1\in H^2 _{l_0} , \label{f1}\\
    v_1 - v_2 &= g_1 \in H^1 (I_2) , \label{g1}\\
    w_1 - w_2 &= h_1\in H^2 _{l_3} , \label{h1}\\
    u_1 ^{(4)}+u_2-\rho _1 u_2'' &=f_2 \in L^2(I_1) , \label{f2}
    \\
    -v_1''+(1+\beta)v_2&=g_2 \in L^2(I_2) ,\label{g2}\\
    w_1 ^{(4)}+w_2-\rho _2 w_2'' &=h_2 \in L^2(I_3) \label{h2} .
\end{align}
Plugging \eqref{f1}--\eqref{h1} in \eqref{f2}--\eqref{h2}, we have to solve
\begin{align}
    u_1 ^{(4)}+u_1-\rho _1 u_1'' &=f_1+f_2-\rho _1 f_1'' ,
    \label{e1}\\
    -v_1''+(1+\beta)v_1&=g_2+(1+\beta )g_1 ,\label{e2}\\
    w_1 ^{(4)}+w_1-\rho _2 w_1'' &=h_1+ h_2-\rho _2 h_1'' \label{e3}.
\end{align}
We define the sesquilinear form $\mathcal{B}:\mathbb{H}\times \mathbb{H} \to \mathbb{C}$ by
\begin{align*}
    \mathcal{B}(Y,\Phi)&:=\langle y_1 , \phi _1 \rangle _{H_{l_0}^2}+\rho _1\langle y' _1 , \phi ' _1 \rangle _{L^2(I_1)}+\langle y_1 , \phi _1 \rangle _{L^2(I_1)}
    +\langle y' _2 , \phi ' _2 \rangle _{L^2(I_2)}\\
    & \quad +(1+\beta)\langle y_2 , \phi _2 \rangle _{L^2(I_2)}
    +\langle y_3 , \phi _3 \rangle _{H_{l_3}^2}+\rho _2\langle y' _3 , \phi ' _3 \rangle _{L^2(I_3)}+\langle y_3 , \phi _3 \rangle _{L^2(I_3)},
\end{align*}
for $Y:=(y_1, y_2, y_3)^\top,\,\Phi:=(\phi _1, \phi _2, \phi _3)^\top \in \mathbb{H}$. 
It is easy to see that $\mathcal{B}:\mathbb{H}\times \mathbb{H} \to \mathbb{C}$ is continuous and coercive.\\

Now, for $(f_1,g_1,h_1,f_2,g_2,h_2)^\top\in \mathcal{H}$, we define  $\Lambda : \mathbb{H} \to \mathbb{C}$ by
\begin{align*}
    \Lambda (\Phi)&:= \langle f_1 +f_2, \phi _1 \rangle_{L^2(I_1)}+\rho _1 \langle f' _1, \phi' _1 \rangle_{L^2(I_1)}+\langle g_2+(1+\beta)g_1 , \phi _2 \rangle_{L^2(I_2)}\\
    &\quad+\langle h_1 +h_2, \phi _3 \rangle_{L^2(I_3)}+\rho _2 \langle h' _1, \phi' _3 \rangle_{L^2(I_3)},
\end{align*}
for all $\Phi \in \mathbb{H}$. It is also easy to see that  $\Lambda : \mathbb{H} \to \mathbb{C}$ is an antilinear continuous functional.\\

Then,  for the Lax-Milgram theorem exists one unique $Y=(y_1,y_2,y_3)^\top\in \mathbb{H}$ such that
\begin{align}
    \mathcal{B}(Y,\Phi)=\Lambda(\Phi),\text{ for all }\Phi \in \mathbb{H}.\label{lm}\end{align}
In particular, if $\phi _1 \in C_{c} ^{\infty} (I_1)$ and $\phi _2 =\phi _3 =0$, we have in (\ref{lm}) that
\begin{align*}
    \langle y_1 , \phi _1 \rangle _{H_{l_0}^2}+\rho _1\langle y' _1 , \phi ' _1 \rangle _{L^2(I_1)}+\langle y_1 , \phi _1 \rangle _{L^2(I_1)}\\
    =\langle f_1 +f_2, \phi _1 \rangle_{L^2(I_1)}+\rho _1 \langle f' _1, \phi' _1 \rangle_{L^2(I_1)},
    \label{y1}
\end{align*}
which can be written, in distibutional sense, as
\begin{align*}
   \langle y ^{(4)} _1 -\rho _1 y'' _1 + y_1 , \phi _1 \rangle & = \langle f_1 +f_2- \rho _1 f'' _1, \phi _1 \rangle,
\end{align*}
for all $\phi _1 \in C_c^\infty(I_1)$. This implies that
$$y ^{(4)} _1 -\rho _1 y'' _1 + y_1=f_1 +f_2- \rho _1 f'' _1 $$
in distributional sense. As $y_1,y_1 '',f_1,f_2, f_1 '' \in L^2(I_1)$, we conclude that $y_1 ^{(4)}\in L^{2}(I_1)$.\\
Analogously, we get  $y_3 ^{(4)}\in L^2(I_3)$ and $y_2 '' \in L^2(I_2)$. \\

Now, let  $(u_1,v_1,w_1):=(y_1,y_2,y_3)$ and $(u_2,v_2,w_2):=(u_1-f_1,v_1-g_1,w_1-h_1)$. Then,  we have that  $U:=(u_1,v_1, w_1,u_2,v_2, w_2)^\top\in \mathbb{H}\times \mathbb{H}$, $u_1 ^{(4)} \in L^2(I_1)$, $v_1 '' \in L^2(I_2)$ y $w_1 ^{(4)} \in L^2(I_3)$. For this $U$ and $\Tilde{U}=(\Tilde{u}_1,\Tilde{v}_1,\Tilde{w}_1,\Tilde{u}_2,\Tilde{v}_2,\Tilde{w}_2)^\top\in \mathbb{H}\times \mathbb{H}$ arbitrary,  
from (\ref{e1})--(\ref{lm}) and $(u_2,v_2,w_2)=(u_1-f_1,v_1-g_1,w_1-h_1)$, we get
\begin{align*}
&\langle -u_1 ^{(4)}+\rho _1 u_2 '' , 
    \Tilde{u}_2 \rangle _{L^2(I_1)}+\langle v_1 '' -\beta v_2 , \Tilde{v}_2 \rangle _{L^2(I_2)}+\langle -w_1 ^{(4)}+\rho _2 w_2 '' ,
    \Tilde{w}_2 \rangle _{L^2(I_3)}\\
    & = \langle u_1 -(f_1+f_2),\Tilde{u}_2\rangle _{L^2(I_1)}+\langle v_1-(g_1+g_2),\Tilde{v}_2 \rangle_{L^2(I_2)}\\
    & \qquad +\langle w_1 -(h_1+h_2),\Tilde{w}_2\rangle _{L^2(I_3)}\\
    & =-\langle u_1 , \Tilde{u}_2 \rangle _{H_{l_0}^2}-\rho _1 \langle (u_1 - f _1)', \Tilde{u}' _2   \rangle_{L^2(I_1)}-\langle v' _1 , \Tilde{v} ' _2 \rangle _{L^2(I_2)}-\beta\langle (v_1 -g_1) , \Tilde{v} _2 \rangle_{L^2(I_2)}\\
    &\qquad-\langle w_1 , \Tilde{w}_2 \rangle _{H_{l_3}^2}-\rho _2 \langle (w_1-h _1)', \Tilde{w}' _2 \rangle_{L^2(I_3)}\\
    &=-\langle u_1 , \Tilde{u}_2 \rangle _{H_{l_0}^2}-\rho _1 \langle u_2 ', \Tilde{u}' _2   \rangle_{L^2(I_1)}-\langle v' _1 , \Tilde{v} ' _2 \rangle _{L^2(I_2)}-\beta\langle v_2 , \Tilde{v} _2 \rangle_{L^2(I_2)}\\
    &\qquad-\langle w_1 , \Tilde{w}_2 \rangle _{H_{l_3}^2}-\rho _2 \langle w_2 ', \Tilde{w}' _2 \rangle_{L^2(I_3)}.
    \end{align*}
Then,
\begin{align*}
    \langle {A}U,\tilde{U}\rangle_{\mathcal{H}} & = \langle u_2 , \Tilde{u}_1 \rangle _{H^2 _{l_0}}+\langle v'_2 , \Tilde{v}'_1 \rangle _{L^2(I_2)}+\langle w_2 , \Tilde{w}_1 \rangle _{H^2 _{l_3}}+\!\langle -u_1 ^{(4)}+\rho _1 u_2 '' , \Tilde{u}_2 \rangle _{L^2(I_1)}\\
  &\quad+\langle v_1 '' -\beta v_2 , \Tilde{v}_2 \rangle _{L^2(I_2)}+\langle -w_1 ^{(4)} + \rho _2 w'' _2 , \Tilde{w}_2 \rangle _{L^2(I_3)}\\
   & = \langle u_2 , \Tilde{u}_1 \rangle _{H^2 _{l_0}}+\langle v'_2 , \Tilde{v}'_1 \rangle _{L^2(I_2)}+\langle w_2 , \Tilde{w}_1 \rangle _{H^2 _{l_3}}-\langle u_1 , \Tilde{u}_2 \rangle _{H_{l_0}^2}\\
   & \quad -\rho _1 \langle u_2 ', \Tilde{u}' _2   \rangle_{L^2(I_1)}-\langle v' _1 , \Tilde{v} ' _2 \rangle _{L^2(I_2)}-\beta\langle v_2 , \Tilde{v} _2 \rangle_{L^2(I_2)}-\langle w_1 , \Tilde{w}_2 \rangle _{H_{l_3}^2}\\
   & \quad -\rho _2 \langle w_2 ', \Tilde{w}' _2 \rangle_{L^2(I_3)}.
\end{align*}
Which means that $U$ satisfies the transmission conditions  (\ref{7})-(\ref{10}) in a weak sense. From this we conclude that $U\in D(\mathcal{A})$, i.e. $b)$ holds. Since $\mathcal{H}$ is a Hilbert space, we get  from $a)$ and $b)$ that  $D(\mathcal{A})$ is dense in $\mathcal{H}$ due to Theorem 4.6, Chapter 1, in \cite{Pazy2012}.
\end{proof}
\begin{theorem}
the operator $\mathcal{A}$ is the generator of a $C_0$-semigroup  $(S(t))_{t\geq 0}$ of contractions on the Hilbert space $\mathcal{H}$. In consequence, for each $U_0 \in D
(\mathcal{A})$ the Cauchy problem  (\ref{cn}) has a unique classical solution $U\in C^1 ([0,\infty),\mathcal{H})$ which depends continuously on the initial data.
\end{theorem}
\begin{proof}
Due to the Proposition \ref{Proposition_wellposedness} and the Lumer-Phillips theorem we have that $\mathcal{A}$ is the generator of a contraction $C_0$-semigroup over $\mathcal{H}$. It follows that, for each $U_0 \in D(\mathcal{A})$, the Cauchy problem (\ref{cn})  has a unique classical solution $U\in C^1 ([0,\infty),\mathcal{H})$ which depends continuously on the initial data, i.e. the problem is well posed.
\end{proof}
\section{Exponential Stability}
\setcounter{equation}{0}
In this section, we will prove the main results of this paper. First we will prove the exponential stability of the semigroup $(S(t))_{t\geq 0}$ generated by $\mathcal{A}$, if we have damping in the three subdomains, i.e., $\rho _1$, $\rho_2$ and $\beta$ are positive.\\
For $U_0 \in D(\mathcal{A})$, we have from Theorem 3.1 that $U(t):=S(t)U_0$ $(t\geq 0)$ is the clasical solution of (\ref{cn}). In this case the energy $E(t)$ of the system  is defined by 
\begin{align*}
    E(t):&=\dfrac{1}{2}||U(t)||^2 _{\mathcal{H}}\\
    &=\dfrac{1}{2}\Big(||u_1 (t)||^2 _{H^2 _{l_0}}+||v _1(t)'||^2_{L^2(I_2)}+||w_1 (t)||^2 _{H^2 _{l_3}}+||u_2 (t)||^2 _{L^2(I_1)}\\
    & \qquad \qquad + ||v _2 (t)||^2_{L^2(I_2)}+ ||w_2 (t)||^2 _{L^2 (I_3)}\Big). 
\end{align*}
Note that
\begin{align}
\dfrac{d}{dt}E(t) & = \text{Re}\langle\mathcal{A} U(t), U(t) \rangle _{\mathcal{H}} \nonumber\\
& = - \rho _1 ||u_2(t)'||^2 _{L^2 (I_1)}-\beta ||v_2  (t)||^2 _{L^2 (I_2)}-\rho _2 ||w_2(t)'||^2 _{L^2 (I_3)},
\end{align}\label{ee}
which shows that the system is dissipative if only one of the damping term is active $(\rho _1 +\rho _2 +\beta >0)$ and conservative if there is no damping at all $(\rho _1 = \rho _2 =\beta = 0)$.\\
Now, we will prove the first main result of this paper.
\begin{theorem}
Let $\rho _1 >0$ , $\rho _2 >0$ and $\beta >0$. Then, the semigroup $(S(t))_{t\geq 0}$ is exponentially stable, i.e., for any $U_0 \in D(\mathcal{A})$ and $U(t):=S(t)U_0$ $(t\geq 0)$ we have
$$E(t)\leq Ce^{-\alpha t}E(0)$$ with positive constants $C$ and $\alpha$.
\end{theorem}
\begin{proof}
For $U_0 \in D(\mathcal{A})$ and $t\geq0$, let $$U(t):=(u_1(t),v_1(t),w_1(t),u_2(t),v_2(t),w_2(t))^{\top}:=S(t)U_0$$ and
$$F(t):=\langle u_1(t) , u_2(t) \rangle_{L^2(I_1)}+\langle v_1(t) , v_2(t) \rangle_{L^2(I_2)}+\langle w_1(t) , w_2(t) \rangle _{L^2(I_3)}.$$ Then,
\begin{align*}
    |F(t)|&\leq \dfrac{1}{2}\Big(||u_1(t)||^2_{L^2 (I_1)}+||u_2 (t)||^2_{L^2 (I_1)}+||v_1 (t)||^2_{L^2 (I_2)}+||v_2 (t)||^2_{L^2 (I_2)}\\
    & \qquad \qquad \qquad + ||w_1 (t)||^2_{L^2 (I_3)} + ||w_2(t)||^2_{L^2 (I_3)}\Big)\\
    &\leq \dfrac{1}{2} ||U(t)||_{X}^2,
\end{align*}
with $X:=H^2(I_1)\times H^1(I_2)\times H^2(I_3) \times \mathbb{L}$. Because of the equivalence between the standar norm in $X$ and the norm $||\cdot||_{\mathcal{H}}$, there exists $c_1 >0$ such that
\begin{align}
    |F(t)| \leq c_1 E(t) .
\end{align}
Due to $\dfrac{d}{dt}U(t)=\mathcal{A}U(t)$, i.e., 
$$\begin{pmatrix}
\dot{u}_1 (t)\\
\dot{v}_1 (t)\\
\dot{w}_1 (t)\\
\dot{u}_2 (t)\\
\dot{v}_2 (t)\\
\dot{w}_2 (t)\\
\end{pmatrix} =
\begin{pmatrix}
u_2 (t)\\
v_2 (t)\\
w_2 (t)\\
-u_1(t)^{(4)}+\rho _1  u_2(t)''\\ 
v_1(t)''-\beta v_2(t)''\\ 
-w_1(t)^{(4)}+\rho _2 w_2(t)''\\ 
\end{pmatrix},
$$
it holds
\begin{align*}
    \dfrac{d}{dt}{F}(t) & =\langle \dot{u}_1 (t) ,u_2 (t) \rangle _{L^2(I_1)}+\langle u_1 (t) ,\dot{u}_2 (t) \rangle _{L^2(I_1)} +\langle \dot{v}_1 (t) ,v_2 (t) \rangle _{L^2(I_2)}\\
    &\qquad +\langle v_1 (t) ,\dot{v}_2 (t) \rangle _{L^2(I_2)} +\langle \dot{w}_1 (t) ,w_2 (t)\rangle _{L^2(I_3)}+\langle w_1 (t) ,\dot{w}_2 (t)\rangle _{L^2(I_3)}\\
         & = ||u_2 (t)||^2 _{L^2(I_1)}+||v_2 (t)||^2 _{L^2(I_2)}+ ||w_2 (t)||^2 _{L^2(I_3)}+\langle \Phi (t) , \mathcal{A}U(t) \rangle _{\mathcal{H}}\\
        & = ||u_2 (t)||^2 _{L^2(I_1)}+||v_2 (t)||^2 _{L^2(I_2)}+ ||w_2 (t)||^2 _{L^2(I_3)}-||u_1 (t)||^2 _{H_{l_0}} \\
        &\qquad -\rho _1 \langle u_1(t)' , u_2(t)'\rangle _{L^2(I_1)}- || v_1(t)'||^2 _{L^2(I_2)}-\beta \langle v_1(t), v_2(t)\rangle _{L^2(I_2)}\\
        & \qquad -||w_1 (t)||^2_{H^2_{l_3}}-\rho _2 \langle w_1(t)', w_2(t)' \rangle _{L^2(I_3)}\\
        &=||(u_2 (t), v_2 (t), w_2 (t))||^2 _{\mathbb{L}}-||(u_1 (t), v_1 (t), w_1 (t))|| ^2 _{\mathbb{H}}\\
        &\qquad-\langle(\rho _1 u_1(t)',\beta v_1 (t), \rho _2 w_1(t)'),(u_2(t)',v_2(t), w_2(t)') \rangle _{\mathbb{L}}.
\end{align*}
There, the weak transmission conditions were used with $\Phi(t):=(0,0,0,u_1 (t),v_1 (t),w_1 (t))^{\top}$.\\
Let $\delta >0$. By Young's inequality and Poincar\'e's inequality, there exists $C_{\delta}>0$  and $c_2>0$ such that
\begin{align*}
-\langle(\rho _1 & u_1(t)',\beta v_1(t), \rho _2 w_1(t)'),(u_2(t)',v_2 (t), w_2(t)') \rangle _{\mathbb{L}}\\
&\leq \delta ||(\rho _1 u_1(t)',\beta v_1(t), \rho _2 w_1 (t)')||^2_{\mathbb{L}}+C_{\delta} ||( u_2(t)', v_2 (t), w_2(t)')||^2_{\mathbb{L}}\\
&\leq c_2\delta ||( u_1  (t), v_1 (t), w_1  (t))||^2_{\mathbb{H}}+C_{\delta} ||( u_2(t)', v_2(t), w_2(t)')||^2_{\mathbb{L}}.\\
\end{align*}
Applying Poincar\'e's inequality to $u_2 (t)$ and $v_2 (t)$, and taking $\delta$ small enough such that $c_2\delta \leq \frac{1}{2}$, we obtain 
\begin{align}
\begin{split}
    \dfrac{d}{dt}F(t)&\leq ||(u_2 (t), v_2 (t), w_2 (t))||^2 _{\mathbb{L}}-\dfrac{1}{2}||(u_1 (t), v_1 (t), w_1 (t))|| ^2 _{\mathbb{H}}\\
    &\quad+ C_{\delta} ||( u_2(t)', v_2 (t), w_2(t)')||^2_{\mathbb{L}}\\
    &\leq c_3||( u_2(t)',v_2 (t), w_2(t)')||^2_{\mathbb{L}}-\dfrac{1}{2}||(u_1 (t), v_1 (t), w_1 (t))|| ^2 _{\mathbb{H}} \label{ff}
\end{split}
\end{align}
for a constant $c_3$. Now, let $L(t):=c_4 E(t)+F(t)$ with $c_4$ a positive constant. Then, for $c_4$ large enough such that $2c_1 \leq c_4$ and $-c_4\text{max}\{\rho _1 ^2 , \beta ^2 , \rho _2 ^2\}+c_3\leq -\frac{1}{2}$, it follows from (\ref{ee}), (\ref{ff}) and Poincar\'e's inequality applyied to $u_2$ and $w_2$ that there are constants $c_5$ and $c_6$ such that
\begin{align}
    \begin{split}
        \dfrac{d}{dt}L(t)&\leq -\dfrac{1}{2}||( u_2(t)',v_2 (t), w_2 (t)')||^2_{\mathbb{L}}-\dfrac{1}{2}||(u_1 (t), v_1 (t), w_1 (t))|| ^2 _{\mathbb{H}}\\
        &\leq -\dfrac{c_5}{2}||( u_2  (t),v_2 (t), w_2  (t))||^2_{\mathbb{L}}-\dfrac{1}{2}||(u_1 (t), v_1 (t), w_1 (t))|| ^2 _{\mathbb{H}}\\
        &\leq -\text{min}\{c_5,1\}\dfrac{1}{2}||U(t)||^2 _{\mathbb{H}\times \mathbb{L}}\\
        &=-c_6 E(t).\label{df}
    \end{split}
\end{align}
As $|F(t)|\leq c_1 E(t)\leq \frac{c_4}{2}E(t)$, we have
$$\dfrac{c_4}{2}E(t)\leq L(t)\leq \dfrac{3c_4}{2}E(t). $$
Therefore, (\ref{df}) yields that $\frac{d}{dt}L(t)\leq -\alpha L(t)$ with some positive constan $\alpha$. By Gronwall's lemma, $L(t)\leq Ce^{-\alpha t}L(0)$, which implies:
$$E(t)\leq \dfrac{2}{c_4}L(t) \leq \dfrac{2}{c_4} e^{-\alpha t}L(0)\leq  e^{-\alpha t}E(0). $$
\end{proof}
Now, we will consider the case in which both beams are undamped. For this we will need to show some results about regularity.
\begin{lemma}\label{Lemma_regularity_1}
Let $I$ an open interval in the real line. For each $g\in L^2(I)$ there exists a unique $v\in H_0^1(I)$ such that 
\begin{equation}\label{Eq_regularity_1}
\int_I v'\,\phi' = - \int_I g\,\phi\qquad (\forall\,\phi\in H_0^1(I)).
\end{equation}
Furthermore, $v\in H^2(I)$.
\end{lemma}
\begin{proof}
It is easy to see that the bilinear form $(v,\phi)\mapsto \int_I v'\,\phi'$ is continuous in $H_0^1(I)\times H_0^1(I)$ and that, due to Poincar\'e's inequality, it is also coercive in $H_0^1(I)$. Then, Lax-Milgram theorem give the existence and uniqueness of the solution $v$ of \eqref{Eq_regularity_1}, since $\phi\mapsto \int_I g\,\phi$ is a continuous linear functional on $H_0^1(I)$, whenever $g\in L^2(I)$. Now, \eqref{Eq_regularity_1} implies that $(v')' = g$ in distributional sense. Since $g\in L^2(I)$ we have that $v'\in H^1(I)$ and therefore $v\in H^2(I)$. 
\end{proof}
We will also need  the result in the following lemma, the proof of which follows analogously to the proof of Corollary 4.3 in \cite{BDHKN19}.
\begin{lemma}\label{Lemma_regularity_2}
Let $a<b$, $f\in L^2((a,b))$ and $z\in\mathbb{C}$. For sufficiently large $\lambda >0$ there exists a unique $u\in H^4((a,b))$ such that
\begin{align*}
u^{(4)} + \lambda u & = f \quad \text{in}\quad (a,b),\\
u(a) = u'(a) & = 0,\\
u''(b) & = 0,\\
u'''(b) & = z.
\end{align*}
\end{lemma}
Now, we can set the following result.
\begin{theorem}
Let $U=(u_1,v_1,w_1,u_2,v_2,w_2)^\top\in D(\mathcal{A})$. Then $u_1\in H^4(I_1)$, $v_1\in H^2(I_2)$ and $w_1\in H^4(I_3)$. In particular, the transmission conditions (\ref{7})--(\ref{10}) hold in the classical sense.
\end{theorem}
\begin{proof}
Let $F=(f_1,g_1,h_1,f_2,g_2,h_2)^\top:=\mathcal{A}U$. Then,
\begin{align}
u_2 & = f_1\in H_{l_0}^2,\label{Eq_Theo_regul_1}\\
v_2 & = g_1\in H^1(I_2),\label{Eq_Theo_regul_2}\\
w_2 & = h_1\in H_{l_3}^2,\label{Eq_Theo_regul_3}\\
-\,u_1^{(4)} + \rho_1u_2'' & = f_2\in L^2(I_1),\label{Eq_Theo_regul_4}\\
v_1'' - \beta v_2 &  = g_2\in L^2(I_2),\label{Eq_Theo_regul_5}\\
-\,w_1^{(4)} + \rho_2 w_2'' & = h_2\in L^2(I_3).\label{Eq_Theo_regul_6}
\end{align}
From \eqref{Eq_Theo_regul_2} and \eqref{Eq_Theo_regul_5} we have 
\begin{equation}\label{Eq_theo_regul_7}
v_1'' = \beta g_1 + g_2.
\end{equation}
Now, let $v_0$ be a smooth function  on $\overline{I_2}$ such that $v_0(l_1)=u_1(l_1)$  and $v_0(l_2)=w_1(l_2)$ (we can take $v_0$ as an affine function for example) and set $\hat{v}:= v_1-v_0$. Then $\hat{v}\in H_0^1(I_2)$ and, due to \eqref{Eq_theo_regul_7}, we have
$$ \hat{v}'' = \beta g_1 + g_2 - v_0''$$
in distributional sense. Then $\hat{v}$ satisfies
$$ \int_{I_2} \hat{v}'\,\phi' = -\,\int_{I_2} (\beta g_1+g_2-v_0'')\phi \qquad (\forall\,\phi\in H_0^1(I_2)).$$
Due to Lemma \ref{Lemma_regularity_1} it follows that $\hat{v}\in H^2(I_2)$ since $\beta g_1+g_2-v_0''\in L^2(I_2)$. Therefore,
$v_1=v_0+\hat{v}\in H^2(I_2)$ and, in particular, Sobolev embedding theorem implies that $v_1\in C^1(\overline{I_2})$.\\
On the other side, we consider $\Phi:=(0,0,0,\phi,\psi,0)^\top$ with $\phi\in H_{l_0}^2$ arbitrary and $\psi\in H^1(I_2)$ such that  $\psi(l_1)=\phi(l_1)$ and $\psi(l_2)=0$. Due to the definition of $D(\mathscr{A})$ we obtain
\begin{align*}
\big\langle \mathscr{A}U , \Phi\big\rangle_{\mathscr{H}} & = \langle -\,u_1^{(4)} + \rho_1 u_2'' , \phi\rangle_{L^2(I_1)} + \langle v_1''-\beta v_2,\psi \rangle_{L^2(I_2)}\\
& = -\,\langle u_1'',\phi''\rangle_{L^2(I_1)} - \rho_1\langle u_2',\phi'\rangle_{L^2(I_1)} - \beta\langle v_2,\psi\rangle_{L^2(I_2)} - \langle v_1',\psi'\rangle_{L^2(I_2)}.
\end{align*}
Clearing de term $\langle u_1^{(4)}, \phi\rangle_{L^2(I_1)}$ in the equality above, taking into account \eqref{Eq_Theo_regul_1} and that $v_1\in H^2(I_2)$, we see that integration by parts yields
\begin{align}\label{Eq_theo_regularity_9}
\langle u_1^{(4)}, \phi\rangle_{L^2(I_1)} & = \langle u_1'',\phi''\rangle_{L^2(I_1)} + \big[\rho_1 f_1'(l_1) - v_1'(l_1)\big]\overline{\phi}(l_1).
\end{align}
By Lemma \ref{Lemma_regularity_2}, for sufficiently large $\lambda>0$ there exists a unique $\tilde{u}_1\in H^4(I_1)$ such that
\begin{align*}
\lambda \tilde{u}_1 + \tilde{u}_1^{(4)} & = \lambda u_1 + \rho f_1'' - f_2 \quad \text{in}\  I_1,\\
\tilde{u}_1(l_0) = \tilde{u}_1'(l_0) & = 0,\\
\tilde{u}_1''(l_1) & = 0,\\
\tilde{u}_1'''(l_1) & = \rho_1 f_1'(l_1) - v_1'(l_1).
\end{align*}
Then, for all $\phi\in H^2_{l_0}$,  we have by integration by parts twice and the boundary conditions above, that
\begin{align}\label{Eq_theo_regularity_10}
\langle \lambda \tilde{u}_1 + \tilde{u}_1^{(4)} , \phi\rangle_{L^2(I_1)} & = \lambda\langle \tilde{u}_1, \phi\rangle_{L^2(I_1)} + \langle \tilde{u}_1'', \phi'' \rangle_{L^2(I_1)}\nonumber\\
& \qquad  \qquad + \big[ \rho_1 f_1'(l_1)-v_1'(l_1)\big]\overline{\phi}_2(l_1).
\end{align}
For the same $\lambda$, adding up the therm $\langle \lambda u_1 , \phi \rangle_{L^2(I_1)}$ in \eqref{Eq_theo_regularity_9} we obtain
\begin{align}\label{Eq_theo_regularity_11}
\langle \lambda {u}_1 + {u}_1^{(4)} , \phi\rangle_{L^2(I_1)} & = \lambda\langle {u}_1, \phi\rangle_{L^2(I_1)} + \langle {u}_1'', \phi'' \rangle_{L^2(I_1)}\nonumber\\
& \qquad  \qquad + \big[ \rho_1 f_1'(l_1)-v_1'(l_1)\big]\overline{\phi}_2(l_1).
\end{align}
From \eqref{Eq_Theo_regul_2} and \eqref{Eq_Theo_regul_4} we have
$$ \lambda u_1 + u_1^{(4)} = \lambda u_1 + \rho_1 f_1'' - f_2 = \lambda \tilde{u}_1 + \tilde{u}_1^{4}. $$
Subtracting \eqref{Eq_theo_regularity_11} from \eqref{Eq_theo_regularity_10}, with $\hat{u}_1:=\tilde{u}_1-u_1$, we get
$$0 = \lambda \langle \hat{u}_1 , \phi \rangle_{L^2(I_1)} + \langle \hat{u}_1'', \phi''\rangle_{L^2(I_1)}$$
for all $\phi\in H_{l_0}^2$. Since $\hat{u}_1\in H_{l_0}^2$ we can set $\phi=\hat{u}_1$ in the last equality and we get $\hat{u}_1 = 0$, which implies $u_1=\tilde{u}_1\in H^4(I_1)$. In particular, Sobolev embedding theorem guarantees that $u_1\in C^3(\overline{I_1})$ and therefore the transmission conditions hold in the classical sense. The proof of $w_1\in H^4(I_3)$, and therefore $w_1\in C^3(\overline{I_3})$, is similar.
\end{proof}
Now, we will use the following  frequency domain result, which  gives us a necessary and sufficient condition fot the exponential stability of a $C_0$-semigroup of contractions. For the proof see \cite{G78}, \cite{H85} and \cite{P84}.
\begin{proposition}\label{Proposition_Pruss}
Let $(T(t))_{t\geq 0}$ be a $C_0$-semigroup of contractions in a Hilbert space $H$, generated by an operator ${A}$. Then the semigroup is exponentially stable if and only if
\begin{equation}\label{Eq_Pruss_condition}
i\mathbb{R} \subset \rho ({A})\quad \text{  and  }\quad ||(i\lambda I -{A})^{-1}||_{\mathcal{L}({H})}\leq C \qquad \forall \lambda \in \mathbb{R}.
\end{equation}
\end{proposition}
The second main result of this paper is the following.
\begin{theorem}
If $\rho_1 = \rho_2 =0$ and $\beta >0$, then the semigroup $(S(t))_{t\geq 0}$ generated by $\mathcal{A}$ is exponentially stable. 
\end{theorem}
\begin{proof}
By the Proposition \ref{Proposition_Pruss} it is sufficient to show that $\mathcal{A}$ satisfies \eqref{Eq_Pruss_condition}. First, we will show that $0 \in  \rho (\mathcal{A}).$ 
Let $F=(f_1, g_1 , h_1 , f_2, g_2 , h_2)^{\top}\in \mathcal{H}$, then we have to find a $U=(u_1, v_1 , w_1 , u_2, v_2 , w_2)^{\top}\in D(\mathcal{A})$ such that $-\mathcal{A}U=F$,
which is equivalent to equations \eqref{Eq_Theo_regul_1}--\eqref{Eq_Theo_regul_6} with $-F$ instead of $F$. Then, plugging the first three equations in the last three, we get
\begin{align}
    u_1 ^{(4)}&=f_2-\rho_1 f_1'',
    \label{re1}\\
   -v_1''&=g_2+\beta g_1 ,
    \label{re2}\\
    w_1^{(4)}&=h_2-\rho_2 h_1''.
    \label{re3}
\end{align}
Now, we define the sesquilinear form $\mathcal{B}_0:\mathbb{H}\times \mathbb{H} \to \mathbb{C}$ by
\begin{align*}
    \mathcal{B}_0 (Y,\Phi)&:=\langle y_1 , \phi _1 \rangle _{H_{l_0}^2}+\langle y' _2 , \phi ' _2 \rangle _{L^2(I_2)}+\langle y_3 , \phi _3 \rangle _{H_{l_3}^2},
\end{align*}
for $Y:=(y_1, y_2, y_3)^\top,\, \Phi:=(\phi _1, \phi _2, \phi _3)^\top\in\mathbb{H}$, and  
 the antilinear functional $\Lambda : \mathbb{H} \to \mathbb{C}$ by
\begin{align*}
    \Lambda (\Phi)&:= \langle f_2, \phi _1 \rangle_{L^2(I_1)}+\rho _1 \langle f' _1, \phi' _1 \rangle_{L^2(I_1)}
    +\langle g_2+\beta g_1 , \phi _2 \rangle_{L^2(I_2)} + \langle h_2, \phi _3 \rangle_{L^2(I_3)}\\
    & \qquad + \rho _2 \langle h' _1, \phi' _3 \rangle_{L^2(I_3)},
\end{align*}
for all $\Phi \in \mathbb{H}$. It is easy to see that $\mathcal{B}_0:\mathbb{H}\times \mathbb{H} \to \mathbb{C}$ is continuous and coercive, and that $\Lambda : \mathbb{H} \to \mathbb{C}$ is continuous. Then, by the Lax-Milgram theorem, there exists a unique  $Y:=(y_1,y_2,y_3)^\top\in \mathbb{H}$ such that
\begin{align}
    \mathcal{B}_0(Y,\Phi)=\Lambda(\Phi),\text{ for all }\Phi \in \mathbb{H}.\label{rlm}
    \end{align}
In the same way as in the proof of Proposition \ref{Proposition_wellposedness}, we obtain that $U:=(y_1,y_2,y_3,-f_1,-g_1,-h_1)^\top\in D(\mathcal{A})$ and satisfies $-\mathcal{A}U = F$. On the other hand, if $\Tilde{U}:=(\tilde{u}_1,\tilde{v}_1,\tilde{w}_1,\tilde{u}_2,\tilde{v}_2,\tilde{w}_2)^\top\in D(\mathcal{A})$ solves $-\mathcal{A}\Tilde{U}=F$, then $\mathcal{B}_0((\tilde{u}_1,\tilde{v}_1,\tilde{w}_1)^\top,\Phi)=\Lambda (\Phi)$ holds for all $\Phi\in \mathbb{H}$ due to definition of $D(\mathcal{A})$ and the weak transmission conditions. Therefore $U=\Tilde{U}$ and $\mathcal{A}$ is a bijection. Since $\mathcal{A}$ is the generator of $C_0$-semigroup of contractions by Theorem 3.1, $\mathcal{A}$ is closed and hence $0\in \rho (\mathcal{A})$.\\
By the Sobolev's embedding theorem, we obtain that  $\mathcal{A}^{-1}$ is a compact operator on $\mathcal{H}$, and therefore, the spectrum of $\mathcal{A}$ consists of eigenvalues only. Thus, we have to establish that there are not purely imaginary eigenvalues. Let $\lambda \in \mathbb{R}$, $\lambda \neq 0$, and $U\in D(\mathcal{A})$ with $\mathcal{A}U=i\lambda U$, i.e.
\begin{align}
    i\lambda u_1 &=u_2, \label{s1}\\
    i\lambda v_1 &=v_2, \label{s2}\\
    i\lambda w_1 &=w_2, \label{s3}\\
    i\lambda u_2+u_1 ^{(4)}&=0, \label{s4}\\
    i\lambda v_2 -v_1''+\beta v_2&=0, \label{s5}\\
    i\lambda w_2+w_1 ^{(4)}&=0. \label{s6}
\end{align}
Due to the dissipativity of $\mathcal{A}$, it holds that
\begin{align*}
0&=\text{Re}\langle (i\lambda -\mathcal{A})U,U \rangle
=\text{Re}(i\lambda ||U||_{\mathcal{H}})-\text{Re}\langle \mathcal{A} U, U \rangle
=\beta ||v_2||^{2} _{L^2(I_2)}.
\end{align*}
Then, $v_1=v_2=0$, i.e. $U=(u_1,0,w_1,u_2,0,w_2)^\top$. Since $v_1 \in H^1 (I_2)$ and $n=1$, the Sobolev's inequality implies that $v_1 (x)=0$ for all $x\in I_1=[l_0,l_1]$.
\\
Multiplying (\ref{s1}) by $i\lambda$ and substituting in (\ref{s4}), we get
\begin{align}
    u_1 ^{(4)}-\lambda u_1=0 \quad\text{ in } (l_0,l_1).\label{s9}
\end{align}
Moreover $u_1$ satisfies the boundary conditions
\begin{align}
    u_1 (l_0)=u_1 ' (l_0)=0\quad \text{ and }\quad u_1''(l_1)=u_1'''(l_1)=0. \label{s10}
\end{align}
Let $x:=l_1+(l_0-l_1)\eta, \quad \eta \in[0,1]$. Then,
$$u_1(x)=u_1(l_1+(l_0-l_1)\eta)=:z(\eta) \quad \text{and}\quad \dfrac{d^{k}u_1}{dx^{k}}=\dfrac{1}{(l_0-l_1)^{k}}\,\dfrac{d^{k}z}{d\eta ^{k}}. $$

Therefore, the problem (\ref{s9})--(\ref{s10}) can be transformed as follows: 
\begin{align}
\begin{split}
    &z^{(4)}-a^2z=0 \qquad \text{in } (0,1),\\
    &z(0)=z''(0)=z'''(0)=0,\label{NEWEQ}\\
    &z(1)=z'(1)=0,
\end{split}
\end{align}
where $a:=(l_0-l_1)^2 |\lambda|\neq 0$. The general solution of the ordinary differential equation in \eqref{NEWEQ} is
$$z (\eta)=c_1 \cosh (\sqrt{a}\eta)+c_2\sinh(\sqrt{a}\eta)+c_3\cos (\sqrt{a}\eta)+c_4\sin (\sqrt{a}\eta).$$
Now, we will see that the boundary conditions in (\ref{NEWEQ}) imply that $c_1=c_2=c_3=c_4=0$ and therefore $z\equiv 0$, which leads to $i\mathbb{R}\subset \rho (\mathcal{A})$. In fact, using the boundary condition in (\ref{NEWEQ}), we obtain
$$ c_1+c_3=0, \quad c_1-c_3=0 \quad \text{and}\quad c_2-c_4=0.$$
Then $c_1=c_3=0$ and $c_2=c_4$. Moreover, $z'(1)=0$ implies $c_2(\cosh(\sqrt{a})+\cos(\sqrt{a}))= 0$. 
Since $\cosh(\sqrt{a})+\cos(\sqrt{a})\neq0$, we have $c_2=0$ and therefore  $z\equiv0$, i.e. $u_1\equiv0$. In similar way we obtain $w_1\equiv0$. Thus, $U=0$ and we conclude that  $i\mathbb{R}\subset \rho (\mathcal{A})$.\\
Now, we will show that
\begin{equation}\label{Eq_condition2_Pruss}
    \sup\limits_{\lambda\in\mathbb{R}}||(i\lambda- \mathcal{A})^{-1}||_{
\mathcal{L}(\mathcal{H})} < \infty.
\end{equation}
If \eqref{Eq_condition2_Pruss} is false, there are sequences  
 $(\lambda _n)_{n\in \mathbb{N}} \subset \mathbb{R}$ and $(U_n)_{n\in\mathbb{N}} \subset D(\mathcal{A})$ such that $|\lambda _n| \xrightarrow[n\rightarrow\infty]{}\infty$,  $||U_n||_{\mathcal{H}} =1$ for all $n\in \mathbb{N}$ and 
\begin{align}
    ||(i\lambda_n - \mathcal{A})U_n||_{\mathcal{H}} \rightarrow 0 \qquad (n\rightarrow \infty).\label{ag1}
\end{align}
Let $F_n := (i\lambda_n - \mathcal{A} )U_n$. Due to the standard norm on  $H^2(I_1)\times H^{1}(I_2)\times H^{2}(I_3)\times L^2(I_1)\times L^2(I_2) \times L^2(I_3)$ is equivalent to the norm $|| \cdot||_\mathcal{H}$ on $\mathcal{H}$, (\ref{ag1}) implies
\begin{align}
    i\lambda_n u_{1,n}-u_{2,n} & = f_{1,n} \rightarrow 0 \text{ in } H^2(I_1),\label{ag2} \\
    i\lambda_n v_{1,n} - v_{2,n} & = g_{1,n} \rightarrow 0  \text{ in } H^1  (I_2),\label{ag3}\\
     i\lambda_n w_{1,n} - w_{2,n} & = h_{1,n} \rightarrow 0 \text{ in } H^2(I_3),\label{ag4}\\
    i\lambda_n u_{2,n}+  u_{1,n}^{(4)} & = f_{2,n}\rightarrow 0 \text{ in } L^2(I_1),\label{ag5}\\
    i\lambda_n v_{2,n} - v_{1,n}'' + \beta v_{2,n} & = g_{2,n} \rightarrow 0  \text{ in }L^2(I_2),\label{ag6}\\
    i\lambda_n w_{2,n} + w_{2,n}^{(4)} & = h_{2,n} \rightarrow 0 \text{ in } L^2(I_3)\label{ag7},
  \end{align}
where $ \rho_1 = \rho_2= 0$. Due to the dissipativity of $\mathcal{A}$, it follows
\begin{align*}
   0\leftarrow \text{Re}(\langle (i\lambda _n - \mathcal{A}) U_n,U_n \rangle) &= \text{Re}[i\lambda_n||U_n||^2_{\mathcal{H}} -\langle \mathcal{A} U_n , U_n \rangle_{\mathcal{H}}] = \beta ||v_{2,n}||^2_{L^2(I_2)},
    \end{align*}
i.e. 
\begin{align}
    ||v_{2,n}||_{L^2(I_2)}\rightarrow 0. \label{ag8}
\end{align}
Now, (\ref{ag3}), (\ref{ag6}) and (\ref{ag8}) imply 
 \begin{equation}\label{ag9}
   ||v_{1,n}||^2 _{L^2(I_2)} \rightarrow 0, \quad   |\lambda_n| ||v_{1,n}||_{L^2(I_2)}\rightarrow 0 \quad \text{and}\quad ||\lambda_n^{-1} v_{1,n}''||_{L^2(I_2)} \rightarrow 0.
 \end{equation}
 From this and the Gagliardo-Nirenberg inequality,  it follows that
 \begin{equation}\label{ag10}
   ||v_{1,n}'||_{L^2(I_2)} \leq 
 ||\lambda_n v_{1,n} ||^{1/2}_{L^2(I_2)}  ||\lambda^{-1}_n v_{1,n}''||^{1/2}_{L^2(I_2)} + ||v_{1,n}||_{L^2(I_2)} \rightarrow 0.   
 \end{equation}
Substituting 
 $ v_{2,n} = i\lambda_n v_{1,n} - g_{1,n}$
 in (\ref{ag6}), we have
 \begin{equation}\label{ag11}
  g_{2,n} =   -\,\lambda_n^2v_{1,n} - i\lambda_ng_{1,n}  - v_{1,n}'' + \beta v_{2,n}.
 \end{equation}
 Now, taking $L^2$-product of \eqref{ag11} with $(l_2 - x)v_{1,n}'$, we get
\begin{align*}
\langle  g_{2,n} & \, , \, (l_2 - x)v_{1,n}'\rangle _{L^2 (I_2)}\\
& = \lambda _{n} \overline{\langle v_{1,n}, (l_2 -x)v_{1,n}'\rangle _{L^2 (I_2)}} \, -\, ||\lambda_n v_{1,n}||^2_{L^2 (I_2)} + (l_2 - l_1)|\lambda_n v_{1,n}(l_1)|^2\\
& \quad - i\lambda_n \langle g_{1,n}\, , \, v_{1,n}\rangle _{L^2 (I_2)} 
+ \langle i\,(l_2 - x)g_{1,n}' \, , \, \lambda_n v_{1,n}\rangle _{L^2 (I_2)}\\
& \quad + i\,(l_2 - l_1) \lambda_n g_{1,n} (l_1) \overline{v_{1,n}}(l_1)-\, ||v_{1,n}'||^2_{L^2 (I_2)} + \overline{\langle v_{1,n}'' \, , \, (l_2 - x) v_{1,n}\rangle _{L^2 (I_2)}}\\
& \quad +(l_2 - l_1)|v_{1,n}'(l_1)|^2 + \langle \beta v_{1,n}\, , \, (l_2 - x)v_{1,n}'\rangle  _{L^2 (I_2)}
\end{align*}
or equivalently
\begin{align}
\begin{split}
-||\lambda_n &v_{1,n} ||^2_{L^2 (I_2)} + (l_2 - l_1) |\lambda_n v_{1,n}(l_1)|^2 - ||v_{1,n}'||^2_{L^2 (I_2)} \\
& \qquad + (l_2 - l_1) |v_{1,n}' (l_1)|^2 + 2\text{Re } \{\langle \beta v_{2,n} \, , \, (l_2 - x) v_{1,n}'\rangle _{L^2 (I_2)}\}\\
& = \langle g_{1,n}\, , \, (l_2 - x)v_{1,n}'\rangle _{L^2 (I_2)} +\overline{ \langle \beta v_{2,n}\, , \, (l_2 -x )v_{1,n}'\rangle _{L^2 (I_2)}}  \\
& \qquad -\overline{\lambda^2 _{n} \langle  v_{1,n} \, , \, (l_2 - x)v_{1,n}'\rangle _{L^2 (I_2)}} + \langle i g_{1,n} \, , \, \lambda _{n} v_{1,n}\rangle _{L^2 (I_2)}   \\
& \qquad - \langle i(l_2 - x) g_{1.n}' \, , \, \lambda_nv_{1,n}\rangle _{L^2 (I_2)} - i(l_2 - l_1)\lambda_n g_{1,n} (l_1) \overline{v_{1,n}}(l_1) \\
& \qquad - \overline{\langle v_{1,n}'' \, , \, (l_2 - x)v_{1,n}'\rangle _{L^2 (I_2)}}. \label{ag12}
\end{split}
\end{align}
From (\ref{ag11}) we have that $\beta v_{2,n}=g_{2,n} + \lambda_n^2v_{1,n} + i\lambda_n g_{1,n} + v_{1,n}''$ and therefore
\begin{align*}
\langle \beta & v_{2,n} \, , \, (l_2 - x) v_{1,n}'\rangle _{L^2 (I_2)}\\
 & =\langle g_{2,n} \, , \, (l_2-x)v_{1,n}'\rangle_{L^2 (I_2)} + \lambda^ 2_n \langle v_{1,n}\, , \, (l_2 - x)v_{1,n}'\rangle _{L^2 (I_2)}\\
&\quad + \langle i g_{1,n}'\, , \,\lambda_n v_{1,n}\rangle _{L^2 (I_2)}-\langle i(l_2-x)g_{1,n}'\, , \, \lambda_n v_{1,n}\rangle _{L^2 (I_2)}\\
& \quad -i(l_2-l_1)\lambda_n g_{1,n} (l_1)\overline{v_{1,n}}(l_1) + \langle v_{1,n}''\, , \, (l_2-x)v_{1,n}'\rangle _{L^2 (I_2)}.
\end{align*}
Substituting this in (\ref{ag12}), we obtain
\begin{align*}
(l_2-&l_1)|\lambda_n v_{1,n}(l_1)|^2 + (l_2-l_1) |v_{1,n}'(l_1)|^2 \\
& = 2\text{Re}\Big\{\langle g_{2,n}\, , \, (l_2-x)v_{1,n}'\rangle _{L^2 (I_2)} - \langle i(l_2-x)g_{1,n}'\, , \, \lambda_n v_{1,n}\rangle _{L^2 (I_2)}\\
&\quad -i(l_2-l_1)g_{1,n}(l_1)\lambda_n\overline{v_1}(l_1) - \langle \beta v_{2,n} \, , \, (l_2-x)v_{1,n}'\rangle _{L^2 (I_2)}\Big\} \\
& \quad + \overline{\langle ig_{1,n}'\, , \,\lambda_nv_{1,n}\rangle _{L^2 (I_2)}} + \langle ig_{1,n}\, , \,\lambda_nv_{1,n}\rangle _{L^2 (I_2)} + ||\lambda_n v_{1,n}||^2_{L^2 (I_2)}\\
& \quad  + ||v_{1,n}'||^ 2_{L^2 (I_2)}. 
\end{align*}
Now, by the Cauchy-Schwarz inequality and Young's inequality, for each $ \epsilon > 0$ there exists a $C_\epsilon>0$ such that
\begin{align*}
(l_2-&l_1)|\lambda_n v_{1,n}(l_1)|^2 + (l_2-l_1) |v_{1,n}'(l_1)|^2 \\
& \leq 2\Big[l_2 ||g_{2,n}||_{L^2 (I_2)} ||v_{1,n}||_{L^2 (I_2)} + l_2 ||g_{1,n}'||_{L^2 (I_2)} ||\lambda_n v_{1,n}||_{L^2 (I_2)}+ \\
&\quad +(l_2-l_1)(\varepsilon |\lambda _n v_{1,n}(l_1)|^2+C_{\varepsilon}|g_{1,n}(l_1)|^2) + \beta l_2 \Vert v_{2,n}\Vert_{L^2(I_2)}\Vert v_{1,n}'\Vert_{L^2(I_2)}\Big]\\
& \quad + 2\Vert g_{1,n}\Vert_{H^1(I_2)}\Vert\lambda_n v_{1,n}\Vert_{L^2(I_2)} + ||\lambda_n v_{1,n}||^2_{L^2 (I_2)} + ||v_{1,n}'||^2_{L^2 (I_2)}.
\end{align*}
From this, \eqref{ag3}, \eqref{ag6}, \eqref{ag8}--\eqref{ag10} and the trace theorem, it follows $|\lambda_n v_{1,n}(l_1)| \xrightarrow[n\rightarrow\infty]{}0$ and $ |v_{1,n}' (l_1)| \xrightarrow[n\rightarrow\infty]{}0$. Thus, 
\begin{equation}\label{ag15}
|\lambda_n u_{1,n}(l_1)| \xrightarrow[n\to\infty]{}0 \, \, \,\,\,  \text{and}  \, \, \,\,\,  |u_{1,n}''' (l_1)|\xrightarrow[n\to\infty]{}0 
\end{equation}
due to the transmission conditions \eqref{6} and \eqref{7} with $\rho_1=0$.\\

Now, substituting $u_{2,n} = i\lambda_n u_{1,n}-f_{1,n}$ in  (\ref{ag5}), we obtain
\begin{equation}\label{ag16}
    f_{2,n} = -\lambda_n^2 u_{1,n} - i\lambda_n f_{1,n} + u_{1,n}^{(4)}. 
\end{equation}
Taking $L^2$-product of (\ref{ag16}) with  $(x-l_0)u_{1,n}'$, we get with integration by parts that
\begin{align}
\begin{split}
\langle f_{2,n}\,&  , \, (x-l_0)u_{1,n}'\rangle _{L^2(I_1)}\\ 
& = ||\lambda_n u_{1,n}||^2_{L^2(I_1)}+ \lambda_n^2 \langle \overline{u_{1,n}\, , \,(x-l_0)u_{1,n}'\rangle _{L^2(I_1)}} - (l_1-l_0)|\lambda_n u_{1,n} (l_1)|^2\\
& \quad + i\lambda_n \langle f_{1,n}\, , \, u_{1,n}\rangle _{L^2(I_1)} + \langle i(x-l_0)f_{1,n}'\, , \,\lambda_n u_{1,n}\rangle _{L^2(I_1)}\\ & \quad -i(l_1-l_0)f_{1,n} \lambda_n \overline{u_{1,n}}(l_1) + \Vert u_{1,n}''\Vert_{L^2(I_1)}^2 - \overline{\langle u_{1,n}''\, , \, (x-l_0)u_{1,n}'''\rangle _{L^2(I_1)}}\\
&\quad + (l_1-l_0)u_{1,n}'''(l_1) \overline{u_{1,n}'(l_1)}. \label{ag17}
\end{split}
\end{align}
From (\ref{ag16}) we have
$ \lambda^2 _n u_{1,n} = - f_{2,n} - i \lambda_n f_{1,n} + u_{1,n}^{(4)}$
and therefore
\begin{align}
\begin{split}\label{agt}
&\lambda^2_{n}\langle u_{1,n} \, , \, (x-l_0)u_{1,n}'\rangle _{L^2(I_1)} \\
& = -\langle f_{2,n}\, , \, (x-l_0)u_{1,n}'\rangle _{L^2(I_1)}
+ \langle if_{1,n}\, , \, \lambda_n u_{1,n}\rangle _{L^2(I_1)}\\
& \quad + \langle i(x-l_0) f_{1,n}'\, , \, \lambda_n u_{1,n}\rangle _{L^2(I_1)} - i (l_1 - l_0) f_{1,n} (l_1) \lambda_n \overline{u_{1,n}(l_1)} \\
& \quad +2||u_{1,n}''|| ^2 _{L^2(I_1} + \langle u_{1,n}''\, , \, (x-l_0) u_{1,n}'''\rangle _{L^2(I_1)}+ (l_1-l_0)u_{1,n}'''(l_1)\overline{u_{1,n}'(l_1)}.
\end{split}
\end{align}
Now, plugging (\ref{agt}) in (\ref{ag17}), we get
\begin{align}
\begin{split}
    ||\lambda _n & u_{1,n}||^2_{L^2(I_1)}+3||u_{1,n}''||^2 _{L^2(I_1)}\\
    &=2\text{Re}\big\{\langle f_{2,n}\, ,\, (x-l_0)u_{1,n}' \rangle _{L^2(I_1)} -\langle if_{1,n}\, ,\,\lambda_n u_{1,n}\rangle _{L^2(I_1)} \\
    & \quad - \langle i(x-l_0) f_{1,n}'\, ,\, \lambda_n u_{1,n}\rangle _{L^2(I_1)} + i(l_1-l_0) f_{1,n} (l_1) \lambda_n \overline{u_{1,n}(l_1)}\\
    & \quad -(l_1-l_0)u_{1,n}'''(l_1)\overline{u_{1,n}'(l_1)}\big\}
    +(l_1 - l_0) |\lambda_n u_{1,n} (l_1)|^2. \label{ag18}
\end{split}
\end{align}
Note that the Gagliardo-Nirenberg inequality implies
$$||u_{1,n}'||_{L^2(I_1)}\leq || u_{1,n}|| ^{1/2}_{L^2(I_1)}||u_{1,n}''|| ^{1/2} _{L^2(I_1)}+||u_{1,n}||_{L^2(I_1)} $$
and thus
\begin{equation}
    ||u_{1,n}'||^2_{L^2(I_1)}\leq 3|| u_{1,n}|| ^2 _{L^2(I_1)}+||u_{1,n}''||_{L^2(I_1)}. \label{ag19}
\end{equation}
Moreover, it follows from the trace theorem that there exists a positive constant $C$ such that 
\begin{equation}
    |u_{1,n}'(l_1)|\leq C||u_{1,n}||_{H^2(I_1)}\leq C||U_{n}||_{\mathcal{H}}= C. \label{ag20}
\end{equation}
Let  $\varepsilon _1$, $\varepsilon _2$ and $\varepsilon _3$ positive numbers. By Young's inequality in (\ref{ag18}), there are positive constants $C_{\varepsilon _1}$,   $C_{\varepsilon _2}$ and $C_{\varepsilon _3}$ such that 
\begin{align}
    \begin{split}
        & ||\lambda _n  u_{1,n} ||^2 _{L^2 (I_1)}+ 3||u_{1,n}''||^2 _{L^2 (I_1)}\\
        & \leq \varepsilon _1  ||u_{1,n}' ||^2 _{L^2 (I_1)}+C_{\varepsilon _1}||f_{2,n}||^2 _{L^2(I_1)} + \varepsilon _2  ||\lambda _n u_{1,n} ||^2 _{L^2 (I_1)}+C_{\varepsilon _2}||f_{1,n}||^2 _{L^2(I_1)}\\
        & \quad +\varepsilon _3  ||\lambda _{n} u_{1,n} ||^2 _{L^2 (I_1)}+C_{\varepsilon _3}||f_{1,n}'||^2 _{L^2(I_1)} + 2(l_1-l_0)\big\{|f_{1,n}(l_1)||\lambda _n u_{1,n}(l_1)|\\
        & \quad + |u_{1,n}'''(l_1)||u_{1,n}'(l_1)|+|\lambda _{n}u_{1,n}(l_1)|^2\big\}\\
    &\leq 3\varepsilon _1  ||\lambda _n u_{1,n} ||^2 _{L^2 (I_1)}+\varepsilon _1  ||u_{1,n}'' ||^2 _{L^2 (I_1)}+C_{\varepsilon _1}||f_{2,n}||^2 _{L^2(I_1)}\\
    & \quad +\varepsilon _2  ||\lambda _n u_{1,n} ||^2 _{L^2 (I_1)}+C_{\varepsilon _2}||f_{1,n}||^2 _{L^2(I_1)} +\varepsilon _3  ||\lambda _{n} u_{1,n} ||^2 _{L^2 (I_1)} + C_{\varepsilon _3}||f_{1,n}'||^2 _{L^2(I_1)}\\
    & \quad +2(l_1-l_0)\big\{|f_{1,n}(l_1)||\lambda _n u_{1,n}(l_1)|+C|u_{1,n}'''(l_1)|+|\lambda _{n}u_{1,n}(l_1)|^2\big\}.\label{ag21}
    \end{split}
\end{align}
There we have used \eqref{ag19} and (\ref{ag20}). Choosing $\varepsilon _1, \varepsilon _2$ and $\varepsilon _3$ small enough such that 
$ 3\varepsilon _1 +\varepsilon _2 + \varepsilon_ 3<1/2$, 
we get from \eqref{ag2}, \eqref{ag5}, (\ref{ag15}) and (\ref{ag21}) that
\begin{align}
    ||\lambda _{n} u_{1,n}||_{L^2 (I_1)}\xrightarrow[n\rightarrow\infty]{}0 \, \, \,\,\, \text{and} \, \, \, \,||u_{1,n}''||_{L^2(I_1)}\xrightarrow[n\rightarrow\infty]{}0.\label{ag22}
\end{align}
Analogously we conclude that 
\begin{align}
    ||\lambda _{n} w_{1,n}||_{L^2 (I_3)}\xrightarrow[n\rightarrow\infty]{}0 \, \, \,\,\, \text{and} \, \, \, \,||w_{1,n}''||_{L^2(I_3)}\xrightarrow[n\rightarrow\infty]{}0.\label{ag23}
\end{align}
Therefore, $||U_{n}||_{\mathcal{H}}\xrightarrow[n\rightarrow\infty]{}0$ due to (\ref{ag22}), (\ref{ag23}), (\ref{ag2}), (\ref{ag3}), (\ref{ag8}) and (\ref{ag10}), which is a contradiction. Thus we have proved that the $(S(t))_{t\geq0}$ is exponentially stable.
\end{proof}

\end{document}